\documentclass{article}		 % Specifies the document class

\title{ON BINARY QUADRATIC FORMS AND THE HECKE GROUPS
	  }

\author{Wendell Culp-Ressler \\
	  Department of Mathematics \\
	  Franklin \& Marshall College \\
	  Lancaster, PA 17604-3003}

\date{August 7, 2002}

\usepackage{amsmath,amsthm,amsfonts,amssymb,amsopn}

   \theoremstyle{plain} 				% This is the default
   \newtheorem{theorem}{Theorem}
   \newtheorem*{corollary}{Corollary}
   \newtheorem{lemma}{Lemma}

   \theoremstyle{definition}
   
   \newtheorem{example}{Example}
   
   \theoremstyle{remark}
   \newtheorem*{remark}{Remark}

   \newcommand{\Z}{\mathbb{Z}}	% The integers
   \newcommand{\Q}{\mathbb{Q}}	% The rationals
   \newcommand{\R}{\mathbb{R}}	% The reals
   \newcommand{\C}{\mathbb{C}}	% The complex numbers
   
\begin{document}

\maketitle

%\setcounter{page}{2}

% \pagebreak

%Proposed running head:  \( \mathbb{Z}[\lambda] \)-BINARY QUADRATIC FORMS

%Send proofs to: 

%\hspace{\parindent}       Wendell Culp-Ressler 

%\hspace{\parindent}	  Department of Mathematics 

%\hspace{\parindent}	  Franklin \& Marshall College 

%\hspace{\parindent}	  Lancaster, PA 17604-3003

%\pagebreak

\begin{abstract}
    
We present a theory of reduction of binary 
quadratic forms with coefficients in \( \mathbb{Z}[\lambda] \), 
where \( \lambda \) is the minimal translation in a Hecke group.  
We generalize from the modular group \( \Gamma(1) = \mathrm{SL}(2,\mathbb{Z}) \)
to the Hecke groups and make extensive use of modified negative continued 
fractions.
We also caracterize the ``reduced'' and ``simple'' hyperbolic fixed 
points of the Hecke groups.

\end{abstract}

%\setcounter{page}{4}

% Introduction

The modular group and negative continued fractions play key roles in 
the
development of a theory of reduction of binary quadratic forms with 
integer coefficients 
as given in \cite{Zag81}.
In this paper we generalize from the modular group to the Hecke groups
by replacing
\( \bigl( \begin{smallmatrix} 1&1\\ 0&1 \end{smallmatrix} \bigr) \) 
with 
\( \bigl( \begin{smallmatrix} 1&\lambda\\0&1 \end{smallmatrix} \bigr) 
\)
for certain values of \( \lambda \) between \( 1 \) and \( 2 \).
We use Rosen's negative \( \lambda \)-continued fractions 
\cite{Ros54},
which are associated with the Hecke groups.
We also make use of original work done in this direction by Schmidt 
and Sheingorn \cite{SS95}.
The result is a theory of reduction of indefinite binary quadratic forms with 
coefficients in any one of the rings \( \Z[\lambda] \).

Part of the motivation for this investigation is our wish to uncover  
properties of \( \Z[\lambda] \)-binary quadratic forms which will be useful for 
characterizing rational period functions of automorphic 
integrals on the Hecke groups \cite{CR99pp}.

In Sections \ref{sec.Heckegroups} and \ref{sec.CFs} we define the 
Hecke groups and \( \lambda \)-continued fractions.
In Section \ref{sec.BQFs} we introduce the idea of a hyperbolic 
\( \Z[\lambda] \)-binary quadratic form and describe a unique 
association between primitive indefinite forms and certain real numbers.
In Section \ref{sec.reducednumbers} we characterize the ``reduced'' 
hyperbolic points of the Hecke groups.
Section \ref{sec.reduction} contains the main reduction theorem.
In Section \ref{sec.Simple} we characterize ``simple'' forms 
and show that they may be put into cycles corresponding to equivalence 
classes.

\pagebreak
\section{HECKE GROUPS}
\label{sec.Heckegroups}
\hspace{\parindent}

In this section we give basic properties of the Hecke groups.
We use the properties of some elements to define a family of functions 
which we will use later to produce cycles of ``simple'' binary 
quadratic forms.

\subsection{DEFINITIONS}
\label{subsec.HGDefns}
\hspace{\parindent}

Let \(S = S_{\lambda} 
   = \bigl( \begin{smallmatrix} 1&\lambda\\ 0&1 \end{smallmatrix} \bigr) \),
\(T = \bigl( \begin{smallmatrix} 0&-1\\ 1&0 \end{smallmatrix} \bigr)\),
and \(I = \bigl( \begin{smallmatrix} 1&0\\ 0&1 \end{smallmatrix} \bigr)\),
where \( \lambda \) is a positive real number.
Put \(G(\lambda)=\langle S,T \rangle / \{\pm I\} \subseteq \mathrm{SL}(2,\R) \).
Erich Hecke \cite{Hec36} showed that the only 
values of \(\lambda\) for which \(G(\lambda)\) is discrete 
are 
\[\lambda = \lambda_{p} = 2\cos(\pi/p), \] 
for \(p=3, 4, 5, \dotsc\),
and \( \lambda \geq 2 \).
We will focus on the discrete groups with \( \lambda < 2 \),
\emph{i.e.}, those with \( \lambda = \lambda_{p} \), \( p \geq 3 \).
These groups have come to be known as the \emph{Hecke 
groups}, 
and we will denote them by
\(G_{p}=G(\lambda_{p})\) for \(p\geq3\).
The first several of these Hecke groups are 
\(G_{3}=G(1)=\Gamma(1)\) (the modular group), 
\( G_{4}=G(\sqrt{2}) \), 
\(G_{5}=G\left( \frac{1+\sqrt{5}}{2} \right)\), and
\( G_{6}=G(\sqrt{3}) \).

Fix \( p \geq 3 \) and let \( U = U_{\lambda_{p}} = S_{\lambda_{p}}T \in G_{p} \).
The generators of \( G_{p} \) satisfy the two relations 
\[T^{2}=U^{p}=I.\]
The groups with \(\lambda \geq 2\) have only one relation \(T^{2}=I\).

The entries of elements of \( G_{p} \)  are in 
\( \Z [\lambda_{p}] \), which for each \( p \geq 3 \) is the ring of algebraic 
integers for \( \Q(\lambda_{p}) \).
For \( p>3 \), \( \Q(\lambda_{p}) \) has nontrivial units and may 
have a nontrivial class group.
For example, \( h_{\Q(\lambda_{68})} =2 \).
This is known to be the only class number greater than \( 1 \) for 
\( p \leq 73 \)  \cite{Was97}.

For 
\(M = \bigl( \begin{smallmatrix} a&b\\ c&d \end{smallmatrix} \bigr)
    \in G_{p}\),
we have \(ad-bc=1\), 
so \(G_{p}\) is a subgroup of \(\mathrm{SL}(2,\Z [\lambda_{p}])\).
It is well-known that \(G_{3}=\mathrm{SL}(2,\Z[\lambda_{3}])\) 
(\emph{i.e.},\(\Gamma(1)=\mathrm{SL}(2,\Z)\)), 
however for the other Hecke groups 
\(G_{p} \subsetneqq \mathrm{SL}(2,\Z[\lambda_{p}])\).

An element \(M = \bigl( \begin{smallmatrix} a&b\\ c&d \end{smallmatrix} \bigr)
    \in G_{p}\) is 
\emph{hyperbolic} if \( \vert a+d \vert >2 \),
\emph{parabolic} if \( \vert a+d \vert =2 \),
and \emph{elliptic} if \( \vert a+d \vert <2 \). 
Elements of the Hecke group act on \( \C \) as linear fractional 
transformations.
A complex number \( z \) is a \emph{fixed point} of \( M \in G_{p} \) 
if \( Mz = z \), so
\( M = \bigl( \begin{smallmatrix} a&b\\ c&d \end{smallmatrix} \bigr) \)
fixes
\begin{eqnarray}
    z & = & \frac{a-d \pm \sqrt{(d-a)^{2}+4bc}}{2c} \nonumber \\
      & = & \frac{a-d \pm \sqrt{(a+d)^{2}-4}}{2c}.
    \label{eq:FixedPointFormula}
\end{eqnarray}
From this it is clear that hyperbolic elements of \( G_{p} \) each have 
two distinct real fixed points.
Parabolic elements each have one real fixed point and elliptic 
elements each have non real fixed points which are complex conjugates of each 
other.

Since \( G_{p} \) is discrete,
the \emph{stabilizer} of any complex number
\( z \) in \( G_{p} \),
\( \{M \in G_{p} \mid Mz=z\} \)
is a cyclic subgroup of \( G_{p} \)
\cite[page 15]{Leh66}.
Thus the fixed point sets of any two elements of \( G_{p} \) 
are identical or disjoint,
and all elements of a stabilizer have identical fixed points.
Accordingly, we may designate fixed points as 
\emph{hyperbolic}, \emph{parabolic}, or \emph{elliptic} 
according to whether the matrices fixing them are hyperbolic, 
parabolic, or elliptic, respectively.
We also define the
\emph{Hecke conjugate} of any hyperbolic fixed point 
to be the other fixed point of the elements in its stabilizer.

\subsection{THE MAP \( \Phi_{p} \) }
\label{subsec.MapPhi}
\hspace{\parindent}

We need the following Lemma in order to define a function
we will use in Section 
\ref{sec.Simple}.
The main idea of the proof will also be useful later.
We thank the referee for simplifying the proof.

\begin{lemma}\label{lem.ForPhi}
    Fix \( p \geq 3 \) and let \( U = U_{\lambda_{p}} \).
    Suppose that \( \alpha \in \R \cup \{\infty\} \) and \( 1 \leq i \leq p-1 \).
    Then \( \alpha \in \left[U^{p-i+1}(0) , U^{p-i}(0) \right) \) if and only if
    \( TU^{i}(\alpha) \in [0,\infty) \), and
    \( \alpha = U^{p-i+1}(0) \) if and only if \( TU^{i}(\alpha) = 0 \).
\end{lemma}

\begin{proof}
    A calculation shows that 
    \[ 0 = U^{p}(0) < U^{p-1}(0) < \cdots < U^{2}(0) < U(0) 
        = \infty, \]
    so the intervals make sense, and
    \begin{equation}
        [-\infty,\infty) = [-\infty,U^{p}(0)) \cup [U^{p}(0),U^{p-1}(0)) 
	     \cup \cdots \cup [U^{2}(0),U(0))
        \label{eq:linedecomp}
    \end{equation}
    is a disjoint union of \( p \) half-open intervals.
    \( U \) maps each interval to the previous interval,
    the first to the last, and left endpoints to left endpoints.
    Thus \( \alpha \in \left[U^{p-i+1}(0) , U^{p-i}(0) \right) 
         = \left[U^{-i+1}(0) , U^{-i}(0) \right) \) 
    if and only if 
    \( U^{i}(\alpha) \in [-\infty,0) \),
    which is true if and only if 
    \( TU^{i}(\alpha) \in [0,\infty) \),
    with left endpoints corresponding to left endpoints.
\end{proof}

By the proof of Lemma \ref{lem.ForPhi},
there is a \( (p-1) \) to one function
\( \Phi_{p} \)
mapping \( [0,\infty) \) to \( [0,\infty) \) defined by
\[ 
  \Phi_{p}(x) = 
      \begin{cases}
          TUx, & U^{p}(0) \leq x < U^{p-1}(0) \\
          TU^{2}x, & U^{p-1}(0) \leq x < U^{p-2}(0) \\
	  \vdots \\
	  TU^{p-1}x,   & U^{2}(0) \leq x,
      \end{cases}
    \]
    where \( U = U_{\lambda_{p}} \).
This function is given more explicitly by
\[ 
  \Phi_{p}(x) = 
      \begin{cases}
          \frac{x}{1-\lambda x},   
	              & 0 \leq x < 1/\lambda \\
	  \frac{1-\lambda x}{(\lambda^{2}-1)x-\lambda },   
	              & 1/\lambda \leq x < \lambda/(\lambda^{2}-1) \\
	  \vdots \\
	  x-\lambda,   & \lambda \leq x,
      \end{cases}
    \]
where \( \lambda = \lambda_{p} \).
Given \( \alpha \in [0,\infty) \), 
the exponent \( i \) in 
\( \Phi_{p}(\alpha) = TU^{i}\alpha \) is the unique exponent 
between \( 1 \) and \( p-1 \) for which \( TU^{i}\alpha \in [0,\infty) \).

\section{CONTINUED FRACTIONS}
\label{sec.CFs}
\hspace{\parindent}

In this section we define \( \lambda_{p} \)-continued fractions.
We also derive some properties of these continued fractions 
for later use.

Rosen \cite{Ros54} introduced a class of continued fractions closely 
associated with the Hecke groups.
He expanded any real number into a continued fraction using a 
\emph{nearest} integral multiple of \( \lambda \) algorithm.
We will use a modification of Rosen's continued fractions, in which we 
expand real numbers into continued fractions using a \emph{next} 
integral multiple of \( \lambda \) algorithm.
The result will be \emph{negative} (or \emph{backwards}) continued 
fractions, in which every numerator is \( -1 \).
These continued fractions, also used in \cite{SS95}, 
generalize the simple negative 
continued fractions in 
\cite{Zag81} for which \( \lambda=1 \) and \( p=3 \).
They will help us make connections between 
poles of RPFs, elements of \( G_{p} \), and binary quadratic forms 
on \(\Z[\lambda_{p}]\).

Fix \( p \geq 3 \), 
put \(\lambda = \lambda_{p}\) 
and let \(r_{j}\in \Z\) 
for \(j \geq 0\).
Define a 
\emph{finite \(\lambda_{p} \)-continued fraction} 
(\( \lambda_{p}\)-CF or \( \lambda\)-CF) by
\begin{eqnarray*}
[r_{0};r_{1},\dots,r_{n}] 
& = & r_{0}\lambda -\cfrac{1}{r_{1}\lambda - 
			  \cfrac{1}{\ddots - 
			    \cfrac{1}{r_{n}\lambda
			}}}			\\
& = & \left( S^{r_{0}}TS^{r_{1}}T \dots S^{r_{n}}T \right) (\infty).
\end{eqnarray*}
We define an \emph{infinite \(\lambda_{p}\)}-CF by
\begin{displaymath}
[r_{0};r_{1},\dots] = \lim_{n \rightarrow \infty} 
[r_{0};r_{1},\dots,r_{n}],
\end{displaymath}
if the limit exists.

We expand any finite real number \( \alpha \) into a 
\(\lambda\)-CF according to the 
next integral multiple of \( \lambda \) algorithm.
Let \(\alpha_{0}=\alpha\) and for \(j \geq 0\) define 
\begin{equation}
    r_{j} = \left[ \frac{\alpha_{j}}{\lambda} \right] +1,
    \label{eq:PartialQuotients}
\end{equation}
and the 
\emph{\( j+1^{\textnormal{st}} \) complete quotient}
\begin{equation}
    \alpha_{j+1}=\frac{1}{r_{j}\lambda - \alpha_{j}}.
    \label{eq:CompleteQuotients}
\end{equation}
Here \( \left[ \: \cdot \: \right] \) is the greatest integer function.
Then \( \alpha_{j}=r_{j}\lambda-\frac{1}{\alpha_{j+1}} \) 
for \( j \geq 0 \) 
and \( [r_{0};r_{1},\dots]\)
is the unique \(\lambda\)-CF for \( \alpha \),
while \( [r_{j};r_{j+1},\dots]\)
is the unique \(\lambda\)-CF for \( \alpha_{j} \).
We note that \eqref{eq:PartialQuotients} 
and \eqref{eq:CompleteQuotients} 
imply that for \( j \geq 1 \) we have
\( r_{j} \geq 1 \) and
\( \alpha_{j} \geq \frac{1}{\lambda} \).

We define an \emph{admissible} \(\lambda\)-CF to be one that arises 
from a finite real number by the
next integral multiple of \( \lambda \) algorithm.
Then we have 
\begin{lemma}\label{lem:admissibleconverge}
    Fix \( p \geq 3 \) and put \( \lambda = \lambda_{p} \).
    Then every admissible \(\lambda\)-CF converges.
\end{lemma}
\begin{proof}
    Put \( S = S_{\lambda} \) and
    let \( [r_{0};r_{1},\ldots] \) be the admissible 
    \(\lambda\)-CF for \( \alpha \in \R \).
    For any \( n \geq 0 \) we define 
    \( C_{n} = [r_{0};r_{1},\ldots,r_{n}] \),
    the \( n \)\emph{th convergent} of the \( \lambda \)-CF.
    We will show that \( \left\{C_{n}\right\}_{n=0}^{\infty} \) 
    is decreasing and bounded below by \( \alpha \).
    
    We define 
    \( C_{m,n} = [r_{m};r_{m+1},\ldots,r_{n}] \)
    for \( 0 \leq m \leq n \) and
    note that \( C_{n,n} = r_{n}\lambda > \alpha_{n} \).
    For all \( n>0 \), \( \alpha_{n}>0 \) and 
    \( C_{n-1,n} = S^{r_{n-1}}T\left(C_{n,n}\right) 
                 > S^{r_{n-1}}T\left(\alpha_{n}\right)
         	 = \alpha_{n-1} \),
    since \( S^{j}T(x) = j\lambda-\frac{1}{x} \) increases 
    monotonically on \( (0,+\infty) \).
    Continuing, we have that 
    \( C_{m,n} > \alpha_{m} \) for \( 0 \leq m \leq n \).
    In particular, 
    \( C_{n} = C_{0,n} > \alpha_{0} = \alpha \) for all \( n \geq 0 \).
    
    In order to show that 
    \( \left\{C_{n}\right\}_{n=0}^{\infty} \) 
    is decreasing, we fix \( n \geq 0 \) and note that
    \( C_{n,n} = r_{n}\lambda 
               > r_{n}\lambda - \frac{1}{r_{n+1}\lambda}
               = C_{n,n+1}\).
    Then for \( n > 0 \),  
    \( C_{n-1,n} = S^{r_{n-1}}T\left(C_{n,n}\right) 
                 > S^{r_{n-1}}T\left(C_{n,n+1}\right)
         	 = C_{n-1,n+1} \).
    Continuing, we have that
    \( C_{m,n} > C_{m,n+1} \) for all 
    \( m \), \( 0 \leq m \leq n \).
    In particular,
    \( C_{n} > C_{n+1} \) for all \( n \geq 0 \).
\end{proof}

A \emph{periodic} \( \lambda_{p} \)-CF is an infinite \( \lambda \)-CF 
that repeats, \textit{i.e.}, one for which there exist \( n \geq 0 \) 
and \( m \geq 1 
\) such that \( r_{j+m}=r_{j} \) for all \( j \geq n \).
We will take \( n \) and \( m \) to be the smallest integers for which 
this happens, and write a periodic \( \lambda \)-CF as 
\[ [r_{0};r_{1},\ldots,r_{n-1},\overline{r_{n},\ldots,r_{n+m-1}}]. \]
A \emph{purely periodic} \( \lambda \)-CF has \( n = 0 \).

Periodic \( \lambda \)-CFs identify the fixed points of \( G_{p} \).
Indeed, we have the following result of Schmidt and Sheingorn
\cite{SS95}.
\begin{lemma}\label{lem.PeriodicCF-FixedPt}
    A real number is a fixed point of \( G_{p} \), \( p\geq 3 \), 
    if and only if it has a periodic \( \lambda_{p} \)-CF expansion.
    Moreover, such a number is parabolic if and only if its 
    \( \lambda_{p} \)-CF has the period 
    \( [\overline{2,\underbrace{1,\ldots,1}_{p-3}}] \), and
    is hyperbolic if and only if its 
    \( \lambda_{p} \)-CF has a period other than
    \( [\overline{2,\underbrace{1,\ldots,1}_{p-3}}] \).
\end{lemma}
\begin{proof}
    This is contained in Lemmas 1, 2, and 3 in \cite{SS95}.
\end{proof}

\pagebreak
There are restrictions on admissible \( \lambda \)-CFs.
We have 
\begin{lemma}\label{lem.AdmissibleLambdaCFs}
    Put \( p \geq 3 \).
    An admissible \( \lambda_{p} \)-CF 
    \begin{enumerate}
    \renewcommand{\labelenumi}{(\roman{enumi})}
        \item  has at most \( p-3 \) consecutive ones in any position 
        but the beginning, and
    
        \item  has at most \( p-2 \) consecutive ones at the beginning.
    \end{enumerate}
\end{lemma}

\begin{remark}\label{rem:AdmissibleCFs}
    In the classical case
    (\( p=3 \) and \( \lambda =1 \))
    Lemma \ref{lem.AdmissibleLambdaCFs} 
    reduces to the fact that admissible 
    simple negative continued fractions 
    have no ones in any position.
    This agrees with the theory of simple negative continued 
    fractions as developed in \cite{Zag81}.
\end{remark}

\begin{proof}
    Set \( \lambda = \lambda_{p} \),
    \( S = S_{\lambda} \) and 
    \( U = U_{\lambda} \), 
    and let \( \alpha \in \R \).
    For any \( m \geq 0 \) and \( j \geq 1 \) we have that 
    \( \alpha_{m} = S^{r_{m}}TS^{r_{m+1}}T \cdots S^{r_{m+j-1}}T\alpha_{m+j} \).
    If the 
    \( \lambda \)-CF has \( j \) 
    consecutive ones starting with 
    \( r_{m} \), then
    \( \alpha_{m} = U^{j}\alpha_{m+j} \).
    By the next multiple of 
    \( \lambda \) algorithm we have that for any
    \( n \geq 1 \),
    \( \alpha_{n} \in [1/\lambda,\infty) \), which is the union of the last 
    \( p-2 \) intervals in \eqref{eq:linedecomp}.
    Thus for \( m \geq 1 \) we must have 
    \( j \leq p-3 \).
    Otherwise  
    \( \alpha_{m+k} \notin [1/\lambda,\infty) \) 
    for some \( k \),
    \( 1 \leq k \leq j \),
    since \( U \) maps each interval to the previous interval in 
    \eqref{eq:linedecomp}.
\end{proof}

The restrictions in Lemma \ref{lem.AdmissibleLambdaCFs} 
are the best possible, since
for any 
\( \lambda = \lambda_{p} \),
\( p \geq 3 \),
\( \frac{3}{2}\lambda + \frac{1}{2}\sqrt{\lambda^{2}+4}
     = [\overline{3;\underbrace{1,\ldots,1}_{p-3}}] 
     \)
and
\( \left(ST\right)^{-2}[\overline{3;\underbrace{1,\ldots,1}_{p-3}}]
     = [1,\overline{\underbrace{1,\ldots,1}_{p-3},3}] 
      \).

\section{BINARY QUADRATIC FORMS}
\label{sec.BQFs}
\hspace{\parindent}

We will exploit the connections between binary quadratic 
forms over \( \Z[\lambda_{p}] \), 
elements of the Hecke group \( G_{p} = G(\lambda_{p}) \), 
and \( \lambda_{p} \)-CFs.

We consider binary quadratic forms with coefficients in  
\( \Z[\lambda_{p}] \),
\[ Q(x,y) = Ax^{2}+Bxy+Cy^{2}. \]
We denote such a form by \( Q=[A,B,C] \) and refer to it as a
\( \lambda_{p} \)-BQF.
We restrict our attention to indefinite forms,
which have positive discriminant 
\( D=B^{2}-4AC \).

\subsection{FORMS AND NUMBERS}
\label{subsec.FormsandNumbers}
\hspace{\parindent}

We associate the 
\( \lambda_{p} \)-BQF \( Q = [A,B,C] \) 
with the number
\( \alpha_{Q} = \frac{-B+\sqrt{D}}{2A} \in \Q(\lambda_{p},\sqrt{D}) \), 
one of the zeros of \( Q(x) = Q(x,1) \).
Under this association, the form 
\( -Q = [-A,-B,-C] \) maps to the other zero of 
\( Q(x) \),
\( \alpha_{-Q} = \frac{B+\sqrt{D}}{-2A} = \frac{-B-\sqrt{D}}{2A} \).
This association is ambiguous in the other direction.
However,
if we restrict out attention to hyperbolic fixed points, we may 
describe a unique association of number to forms.
We thank the referee for providing a shorter proof 
of the following Lemma.

\begin{lemma}\label{lemma.NumbertoForm}
    Every hyperbolic fixed point of \( G_{p} \),
    \( p \geq 3 \), may be 
    associated with a unique indefinite 
    \( \Z[\lambda_{p}] \)-binary quadratic form.
\end{lemma}

\begin{proof}
    Suppose that 
    \( M = \bigl( \begin{smallmatrix} a&b\\ c&d \end{smallmatrix} \bigr) \)
    and
    \( M^{-1} = \bigl( \begin{smallmatrix} d&-b\\ -c&a \end{smallmatrix} \bigr) \)
    generate the stabilizer of hyperbolic fixed point 
    \( \alpha \) in \( G_{p} \).
    Now \( M\alpha = \alpha \) gives rise to the equation
    \[ c\alpha^{2}+(d-a)\alpha-b=0, \]
    and \( M^{-1}\alpha = \alpha \) gives rise to the equation
    \[ -c\alpha^{2}-(d-a)\alpha-(-b)=0. \]
    We put
        \begin{displaymath}
	M_{\alpha} = 
	\begin{cases}
	    \bigl( \begin{smallmatrix} a&b\\ c&d \end{smallmatrix} \bigr) & 
	        \text{if } \alpha = \frac{a-d+\sqrt{D}}{2c}, \\
            \bigl( \begin{smallmatrix} d&-b\\ -c&a \end{smallmatrix} \bigr) & 
	        \text{if } \alpha = \frac{a-d-\sqrt{D}}{2c},
	\end{cases}
        \end{displaymath}
    and
        \begin{displaymath}
	    \label{eq:UniqueBQF1}
            Q_{\alpha} = 
	    \begin{cases}
	        [c,d-a,-b] & 
	            \text{if } \alpha = \frac{a-d+\sqrt{D}}{2c}, \\
	        -[c,d-a,-b] & 
	            \text{if } \alpha = \frac{a-d-\sqrt{D}}{2c},
	    \end{cases}
	\end{displaymath}
    where 
    \( D=(d-a)^{2}+4bc=(a+d)^{2}-4 \).
\end{proof}
If a \( \lambda \)-BQF \( Q \) arises from a hyperbolic fixed point 
as in the proof of
Lemma \ref{lemma.NumbertoForm}, we say that \( Q \) is 
\emph{hyperbolic}.

It is easy to calculate the generators of the stabilizer of 
any fixed point
\( \alpha \) in \( G_{p} \).
By Lemma \ref{lem.PeriodicCF-FixedPt}, 
\( \alpha \) has a periodic \( \lambda_{p} \)-CF expansion
\( \alpha = [r_{0};r_{1},\ldots,r_{n-1},\overline{r_{n},\ldots,r_{n+m-1}}] \).
We put \( V=S^{r_{0}}TS^{r_{1}}T \cdots S^{r_{n-1}}T \) and
\( W=S^{r_{n}}TS^{r_{n+1}}T \cdots S^{r_{n+m-1}}T \).
Then \( M=VWV^{-1} \) and \( M^{-1}=VW^{-1}V^{-1} \) generate the stabilizer of 
\( \alpha \) in \( G_{p} \) since every 
\( \lambda_{p} \)-CF period is minimal.

Suppose that \( \alpha = \frac{-B+\sqrt{D}}{2A} \) 
is a hyperbolic fixed point of 
\( G_{p} \) associated with the hyperbolic \( \lambda_{p} \)-BQF
\( Q_{\alpha} = [A,B,C] \).
We denote the Hecke conjugate of 
\( \alpha \)
by
\( \alpha^{\prime} \)
and observe that by the proof of
Lemma \ref{lemma.NumbertoForm},
\( \alpha^{\prime} = \frac{-B-\sqrt{D}}{2A} \) is
associated with 
\( Q_{\alpha^{\prime}}=-Q_{\alpha} = [-A,-B,-C] \).
A straightforward calculation shows that if \( V \in G_{p} \),
then \( \left(V\alpha\right)^{\prime} = V\alpha^{\prime} \).

\subsection{EQUIVALENCE CLASSES}
\label{subsec.EquivalenceClasses}
\hspace{\parindent}

Elements of a Hecke group act on \( \lambda \)-BQFs by
\( \left(Q \circ M \right)(x,y) = Q(a x + b y, c x + d y) \)
for 
\( M = \bigl( \begin{smallmatrix} a &b \\ c & d 
              \end{smallmatrix} \bigr)\ \in G_{p} \).
More explicitly, if \( Q = [A,B,C] \) we have that
\[ [A,B,C] \circ \bigl( \begin{smallmatrix} a &b \\ c & d \end{smallmatrix} \bigr)\
    = [A^{\prime},B^{\prime},C^{\prime}], \]
where 
\[ A^{\prime} = Aa^{2}+Bac+Cc^{2} = Q(a,c), \]
\[ B^{\prime} = 2Aab+B(ad+bc)+2Ccd, \] 
and
\[ C^{\prime} = Ab^{2}+Bbd+Cd^{2} = Q(b,d). \]
A straightforward calculation shows that 
\( B^{\prime 2}-4A^{\prime}C^{\prime} = B^{2}-4AC \),
so the action of a Hecke group preserves the discriminant.

We say that \( Q \) and \( Q^{\prime} \) are
\emph{\( G_{p} \)-equivalent}, and write \( Q \sim Q^{\prime} \),
if there exists a \( V \in G_{p} \) such that 
\( Q^{\prime} = Q \circ V \).
It is easy to check that \( G_{p} \)-equivalence is an equivalence relation, so 
\( G_{p} \) partitions the \( \lambda \)-BQFs into equivalence 
classes of forms.

\begin{lemma}\label{lem:FormsandNumbers}
    Fix \( p \geq 3 \) and let 
    \( \lambda = \lambda_{p} \).
    Suppose that 
    \( Q_{\alpha} \) and \( Q_{\beta} \) are hyperbolic 
    \( \lambda \)-BQFs associated with hyperbolic numbers
    \( \alpha \) and \( \beta \), respectively,
    and let \( V \in G_{p} \).
    Then \( Q_{\beta} = Q_{\alpha} \circ V \) if and only if 
    \( \beta = V^{-1}\alpha \).
\end{lemma}
\begin{proof}
    We write 
    \( Q_{\alpha} = [A,B,C] \), so 
    \( \alpha = \frac{-B+\sqrt{D}}{2A} \), where 
    \( D=B^{2}-4AC \), 
    and we write 
    \( V = \bigl( \begin{smallmatrix} a & b \\ 
                                      c & d 
                  \end{smallmatrix} \bigr)\ \).
		  
    Suppose that \( Q_{\beta} = Q_{\alpha} \circ V \). 
    Then 
    \begin{eqnarray*}
        Q_{\beta} & = & [A,B,C] \circ 
                           \bigl( \begin{smallmatrix} a & b \\ 
                                                    c & d  \\
                                  \end{smallmatrix} 
			   \bigr)\                        \\
         & = & [Aa^{2}+Bac+Cc^{2},
	        2Aab+B(ad+bc)+2Ccd,
		Ab^{2}+Bbd+Cc^{2}],
    \end{eqnarray*}
    so
    \[ \beta = 
        \frac{-(2Aab+B(ad+bc)+2Ccd)+\sqrt{D}}
	    {2(Aa^{2}+Bac+Cc^{2})}. \]
    On the other hand, 
    \( V^{-1} = \bigl( \begin{smallmatrix} d & -b \\ 
                                      -c & a 
                  \end{smallmatrix} \bigr)\ \),
    so 
    \begin{eqnarray*}
        V^{-1}\alpha & = & \frac{d 
        \left(\frac{-B+\sqrt{D}}{2A}\right)-b}
	                        {-c\left(\frac{-B+\sqrt{D}}{2A}\right)+a}  \\
         & = & \frac{-2Ab -Bd +c\sqrt{D}}
	            {2Aa -Bc -c\sqrt{D}}  \\
         & = & \frac{-(2Aab+B(ad+bc)+2Ccd)
	               +\sqrt{D}}
	    {2(Aa^{2}+Bac+Cc^{2})}  \\
	 & = & \beta.
    \end{eqnarray*}
    
    Now suppose that \( \beta = V^{-1}\alpha \).
    Let \( M \) be one of the two generators of the 
    stabilizer of \( \alpha \) in \( G_{p} \).
    Then \( W=V^{-1}MV \) is one of the two generators of the 
    stabilizer of 
    \( \beta \) in \( G_{p} \).
    Moreover, \( M \) and \( W \) have the same trace, 
    since the trace of a matrix is preserved by conjugation.
    Now \( M\alpha = \alpha \) gives rise to \( Q_{\alpha} \) and 
    \( W\beta = \beta \) gives rise to \( Q_{\beta} \) 
    as in the proof of
    Lemma \ref{lemma.NumbertoForm}, so
    \( Q_{\alpha} \) and \( Q_{\beta} \) have the same discriminant.
    By the calculation above, 
    \( Q_{\alpha} \circ V \) has the associated hyperbolic number
    \[ \frac{-(2Aab+B(ad+bc)+2Ccd)+\sqrt{D}}
	    {2(Aa^{2}+Bac+Cc^{2})}
	= V^{-1}\alpha = \beta. \]
    A similar calculation shows that 
    \( -(Q_{\alpha} \circ V) \) has the associated hyperbolic number
    \[ \frac{-(2Aab+B(ad+bc)+2Ccd)-\sqrt{D}}
	    {2(Aa^{2}+Bac+Cc^{2})}
	= V^{-1}\alpha^{\prime} = \beta^{\prime}. \]    
    Thus \( (Q_{\alpha} \circ V)(x,1) \) and \( Q_{\beta}(x,1) \) 
    have the same discriminant and the same zeros, 
    so they are equal.
    Hence \( Q_{\alpha} \circ V = Q_{\beta} \).
\end{proof}

Lemma \ref{lem:FormsandNumbers} implies that
\( G_{p} \)-equivalence of hyperbolic \( \lambda \)-BQFs induces a corresponding 
\( G_{p} \)-equivalence of associated numbers.
This equivalence 
yields the following.

\begin{lemma}\label{lem.EquivFixedPts}
    Two hyperbolic fixed points of 
    \( G_{p} \), \( p \geq 3 \) , are equivalent if and only 
    if the periods of their 
    \( \lambda_{p} \)-CF expansions are cyclic 
    permutations of each other.
\end{lemma}

\begin{corollary}\label{cor.HyperbolicClasses}
    Let \( p \geq 3 \) and put \( \lambda = \lambda_{p} \).
    Every \( G_{p} \)-equivalence class of 
    \( \lambda \)-BQFs contains either all 
    hyperbolic forms or no hyperbolic forms.
\end{corollary}
\begin{proof}
    \( G_{p} \)-equivalence classes of hyperbolic \( \lambda_{p}\)-BQFs are in 
    \( 1-1 \) correspondence with \( G_{p} \)-equivalence classes of 
    hyperbolic fixed points of \( G_{p} \).
    The Corollary follows from Lemma \ref{lem.PeriodicCF-FixedPt} and
    Lemma \ref{lem.EquivFixedPts}.
\end{proof}
We may now label equivalence classes themselves as hyperbolic or 
non-hyperbolic.

\section{REDUCED NUMBERS}
\label{sec.reducednumbers}
\hspace{\parindent}

In this section we generalize the reduction theory of indefinite binary 
quadratic forms to \(\lambda \)-BQFs.
We characterize ``reduced'' numbers and we 
see that every hyperbolic \( G_{p} \)-equivalence class of forms 
contains a cycle of ``reduced'' forms.
From this it follows that 
there are a finite number of hyperbolic equivalence classes for each 
discriminant.

Fix \( p \geq 3 \).
We say that a real number \( \alpha \) is a \emph{\( G_{p} \)-reduced} 
number if the \(\lambda_{p} \)-CF expansion of \( \alpha \) is purely 
periodic, with period other than 
\( [\overline{2,\underbrace{1,\ldots,1}_{p-3}}] \).
If \( \alpha \) is \( G_{p} \)-reduced,
we also say that the associated \(\lambda_{p} \)-BQF \( Q_{\alpha} \) is
\( G_{p} \)-reduced.

The following Lemma is a modification of a familiar result from classical 
continued fractions.
It is also stated in \cite[page 389]{SS95}.
\begin{lemma}\label{lem.CFofConjugate}
    Let \( p \geq 3 \) and put \( \lambda = \lambda_{p} \).
    If 
    \( \alpha = [\overline{r_{0};r_{1},\ldots,r_{n}}] \)
    is a \(\lambda \)-CF,
    then 
    \( \frac{1}{\alpha^{\prime}} = [\overline{r_{n};r_{n-1},\ldots,r_{0}}]  \).
\end{lemma}
\begin{proof}
    Put \( S = S_{\lambda} \).
    We have \(  \alpha_{j} = [\overline{r_{j};r_{j+1},\ldots}] \) 
    for \( j \geq 0 \).
    Then \( \alpha_{j+1} = TS^{-r_{j}}\alpha_{j} \), 
    and \( \alpha_{0} = \alpha_{n+1} = TS^{-r_{n}}\alpha_{n} \).
    Taking Hecke conjugates, we have 
    \( \alpha_{j+1}^{\prime} = TS^{-r_{j}}\alpha_{j}^{\prime}
                        = \frac{1}{r_{j}\lambda - \alpha_{j}^{\prime}}\), 
    \( j\geq 0 \), and 
    \( \alpha_{0}^{\prime} = \frac{1}{r_{n}\lambda - \alpha_{n}^{\prime}} \),
    so 
    \( \frac{1}{\alpha_{0}^{\prime}} 
        = r_{n}\lambda - \alpha_{n}^{\prime} \).
    We combine these to get
    \begin{eqnarray*}
        \frac{1}{\alpha_{0}^{\prime}}  & 
	    = & r_{n}\lambda - \cfrac{1}{r_{n-1}\lambda-
                  \cfrac{1}{r_{n-2}\lambda-\ddots }}  \\
         & = & [\overline{r_{n};r_{n-1},\ldots,r_{0}}].
    \end{eqnarray*}
    \end{proof}

We need the following result in the proof of the 
next Theorem.
\begin{lemma}\label{lem.Ureciprocals}
    Fix \( p \geq 3 \) and let \( U = U_{\lambda_{p}} \).
    Then \( \frac{1}{U^{k}(0)}=U^{p-k+1}(0) \) 
    for any integer \( k \).
\end{lemma}

\begin{proof}
    Let 
    \( c_{k} = \frac{sin\left(k\pi/p\right)}{sin\left(\pi/p\right)} \) 
    for \( k\geq 0 \).
    Meier and Rosenberger \cite{MR84} show that as linear fractional 
    transformations
    \begin{equation}
        U^{k} 
        = \bigl( \begin{smallmatrix} c_{k+1}&-c_{k}\\ c_{k}&-c_{k-1} \end{smallmatrix}
	  \bigr),
        \label{eq:powerofU}
    \end{equation}
    for \( k \in \Z^{+} \).
    In fact, it is easy to show that \eqref{eq:powerofU}
    holds for all integers \( k \).
    Then 
    \begin{eqnarray*}
        U^{p-k+1}(0) & = & \frac{c_{p-k+1}}{c_{p-k}}  \\
         & = & \frac{sin\left((p-k+1)\pi/p\right)}{sin\left((p-k)\pi/p\right)} \\
         & = & \frac{sin\left((k-1)\pi/p\right)}{sin\left(k\pi/p\right)} \\
         & = & \frac{c_{k-1}}{c_{k}} \\
         & = & \frac{1}{U^{k}(0)}.
    \end{eqnarray*}
\end{proof}

We next characterize \( G_{p} \)-reduced numbers.
We thank the referee for making suggestions which helped to correctly 
formulate and prove the following Theorem.
\begin{theorem}\label{thm.ReducedNumbers}
    Put \( p \geq 3 \),
    \( \lambda = \lambda_{p} \), and
    \( U = U_{\lambda} \).
    Suppose that \( \alpha \) is a hyperbolic fixed point of
    \( G_{p} \).
    Then \( \alpha \) is \( G_{p} \)-reduced
    with \( k \) leading ones in its \( \lambda \)-CF
    if and only if 
    \( k \) is the smallest nonnegative integer such that
    \begin{equation}
        0 < \alpha^{\prime} < U^{k+2}(0) < \alpha < U^{k+1}(0).
        \label{eq:ReducedNumbers}
    \end{equation}
\end{theorem}

\begin{remark}\label{rem:ClassicalReducedNumbers}
    In the classical case 
    (\( p=3 \) and \( \lambda =1 \))
    we must have that \( k=0 \).
    Then a hyperbolic point \( \alpha \) is
    \( \left(\Gamma(1)-\right) \) reduced 
    if and only if
    \( 0<\alpha^{\prime}<1<\alpha<\infty \).
    This agrees with the definition of a reduced number 
    in \cite[page 128]{Zag81}.
\end{remark}

\begin{proof}
    First suppose that \( \alpha \) is \( G_{p} \)-reduced
    and that the 
    \( \lambda \)-CF for \( \alpha \) 
    has \( k \) leading ones.
    Let \( m \in \Z^{+} \) denote the period length of the
    \( \lambda \)-CF for \( \alpha \).
    Then 
    \( k \leq \min\{m-1,p-3\} \).
    We may write
    \( \alpha = [\overline{r_{0};r_{1},\ldots,r_{m-1}}] 
         = [\overline{\underbrace{1;1,\ldots,1}_{k},
	   r_{k},r_{k+1},\ldots,r_{m-1}}] \),
    and
    \( \alpha_{j} = [\overline{r_{j};r_{j+1},\ldots,r_{m-1},
                r_{0},\ldots,r_{j-1}}] \)
    for \( 0 \leq j \leq m-1 \).
    Each \( \alpha_{j} \) 
    is hyperbolic and thus not a multiple of \( \lambda \).
    Since \( r_{k} \geq 2 \) we have
    \( \alpha_{k} > \lambda = U^{2}(0) \),
    \emph{i.e.},
    \( \alpha_{k} \) is in the \( p \)th interval of 
    \eqref{eq:linedecomp}.
    Now \( r_{j}=1 \) for \( 0 \leq j < k \) 
    implies that 
    \( \alpha = U^{k}\alpha_{k} \), so
    \( \alpha \) is in the \( (p-k) \)th interval of 
    \eqref{eq:linedecomp},
    \emph{i.e.},
    \( U^{k+2}(0) < \alpha < U^{k+1}(0) \).
    By Lemma \ref{lem.CFofConjugate}
    we have
    \( \frac{1}{\alpha^{\prime}} 
        = [\overline{r_{m-1};\ldots,r_{k+1},r_{k},
  	   \underbrace{1,\ldots,1}_{k}}]
	> 0 \),
    so \( \alpha^{\prime} >0 \).
    Also, \( \alpha_{k} 
        = [\overline{r_{k};r_{k+1},\ldots,r_{m-1},
	    \underbrace{1,\ldots,1}_{k}}] \),
    so \( \frac{1}{\alpha_{k}^{\prime}} 
        = [\overline{\underbrace{1;1,\ldots,1}_{k},
	    r_{m-1},\ldots,r_{k+1},r_{k}}] 
	> \frac{1}{\lambda} = U^{p-1}(0) \).
    Then \( \frac{1}{\alpha^{\prime}} 
        = U^{-k}\left(\frac{1}{\alpha_{k}^{\prime}}\right) 
	> U^{p-k-1}(0) 
	> 0 \),
    since \( U^{-k} \) maps every interval in
    \eqref{eq:linedecomp} 
    \( k \) intervals to the right, and 
    \( 2 \leq p-k-1 \leq p-1 \).
    Thus
    \( \alpha^{\prime} < \frac{1}{U^{p-k-1}(0)}
                     = U^{k+2}(0)\),
    by Lemma \ref{lem.Ureciprocals},
    and we have verified \eqref{eq:ReducedNumbers}.
    
    Next, suppose that 
    \( k \) is the smallest nonnegative integer such that
    \eqref{eq:ReducedNumbers} holds.
    Since \( U^{p}=I \),
    we have that 
    \( k \leq p-1 \).
    Furthermore
    \( k \leq p-2 \),
    since \( k=p-1 \) implies that 
    \( \infty = U^{p+1}(0) < U^{p}(0) =0 \).
    
    By \eqref{eq:ReducedNumbers}, 
    \( \alpha \) is in the
    \( p-k^{\textnormal{th}} \) interval from the left
    and \( \alpha^{\prime} \) is in one of the second through the
    \( p-k-1^{\textnormal{st}} \) intervals of \eqref{eq:linedecomp}.
    Let \( \alpha = 
    [r_{0};r_{1},\ldots,r_{n-1},\overline{r_{n},\ldots,r_{n+m-1}}] \)
    denote the \( \lambda- \)CF for \( \alpha \).
    A complete quotient
    \( \alpha_{j} \) is in an interval other than the 
    \( p^{\textnormal{th}} \) interval when
    \( 0<\alpha_{j}<\lambda \) and
    \( r_{j}=1 \).
    In this case
    \( \alpha_{j+1} = U^{-1}\alpha_{j} \) and
    \( \alpha_{j+1}^{\prime} = U^{-1}\alpha_{j}^{\prime} \), 
    so 
    \( \alpha_{j+1} \) is in the interval to the right of \( \alpha_{j} \)
    and
    \( \alpha_{j+1}^{\prime} \) is in the interval to the right of 
    \( \alpha_{j}^{\prime} \).

    Since \( \alpha \) is in the
    \( p-k^{\textnormal{th}} \) interval from the left,
    the discussion above implies that
    each \( \alpha_{j} \), \( 0\leq j \leq k-1 \),
    is in an interval to the left of the \( p^{\textnormal{th}} \) interval,
    while 
    \( \alpha_{k} \)
    is in the \( p^{\textnormal{th}} \) interval.
    Furthermore
    \( r_{j}=1 \) for \( 0 \leq j \leq k-1 \)
    and \( r_{k} \geq 2 \),
    \emph{i.e.}, the \( \lambda \)-CF for 
    \( \alpha \) has \( k \) leading ones.

    In order to show that \( \alpha \) is 
    \( G_{p} \)-reduced, 
    we first claim that \( 0 < \alpha_{j}^{\prime} < \lambda \)
    for all \( j \geq 0 \).
    Since \( \alpha^{\prime} \) is in one of the second through the 
    \( p-k-1^{\textnormal{st}} \) intervals of
    \eqref{eq:linedecomp},
    and the \( \lambda \)-CF for \( \alpha \) has \( k \) leading ones,
    the discussion above implies that each 
    \( \alpha_{j}^{\prime} \),
    \( 0 \leq j \leq k \) is in one of the second through the 
    \( p-1^{\textnormal{st}} \) intervals,
    \emph{i.e.}, 
    \( 0 < \alpha_{j}^{\prime} < \lambda \) for 
    \( 0 \leq j \leq k \).
    A calculation shows that if 
    \( 0 < \alpha_{t}^{\prime} < \lambda \) and \( r_{t} \geq 2 \)
    for any \( t \), then
    \( 0 < \alpha_{t+1}^{\prime} < \frac{1}{\lambda} \) and
    \( \alpha_{t+1}^{\prime} \) is in the second interval of
    \eqref{eq:linedecomp}.
    If \( r_{t} \) is followed by 
    \( r_{t+1} \geq 2 \), then 
    \( 0 < \alpha_{t+2}^{\prime} < \frac{1}{\lambda} \)
    by the same calculation.
    On the other hand, 
    if \( r_{t} \) is followed by 
    \( \ell \) ones,
    \( \alpha_{j}^{\prime} \) is in one of the 
    \( 2^{\textnormal{nd}} \) through the 
    \( \ell + 2^{\textnormal{nd}} \) intervals for
    \( t+1 \leq j \leq t + \ell +1 \).
    Since the \( \lambda \)-CF for \( \alpha_{j}^{\prime} \)
    is admissible,
    \( \ell \leq p-3 \), so
    each such
    \( \alpha_{j}^{\prime} \) is in one of the 
    \( 2^{\textnormal{nd}} \) through the 
    \( p-1^{\textnormal{st}} \) intervals.
    Thus \( 0 < \alpha_{j}^{\prime} < \lambda \) for the
    \( \ell \) complete quotients following 
    \( \alpha_{t} \).
    The claim follows by induction on \( j \).

    We next note that for every 
    \( j \geq 1 \),
    the complete quotient
    \( \alpha_{j} \) 
    has a unique predecessor
    \( \alpha_{j-1} \) with
    \( 0 < \alpha_{j}^{\prime} < \lambda \).
    Indeed, 
    \( 0 < \alpha_{j}^{\prime} < \lambda \) implies that 
    \( T\alpha_{j}^{\prime} < 0 \), so
    \( 0 < S^{t}T\alpha_{j}^{\prime} < \lambda \)
    for a unique \( t \), 
    \( t \geq 1 \).
    But since 
    \( 0 < \alpha_{j-1}^{\prime} < \lambda \) and
    \( \alpha_{j-1}^{\prime} = S^{r_{j-1}}T\alpha_{j}^{\prime} \),
    we must have that 
    \( r_{j-1} = t \) is uniquely determined.
    Then 
    \( \alpha_{n} = \alpha_{n+m} \) implies that 
    \( \alpha_{j} = \alpha_{j+m} \) for every 
    \( j < n \)
    and 
    \( \alpha = 
    [\overline{r_{0};r_{1},\ldots,r_{m-1}}] \).
\end{proof}

\section{REDUCTION}
\label{sec.reduction}
\hspace{\parindent}

We are ready to describe the reduction of hyperbolic elements of \( G_{p} \).
\begin{theorem}\label{thm.ReductionofNumbers}
    Let \( p \geq 3 \) and put
    \( \lambda = \lambda_{p} \).
    Every hyperbolic fixed point \( \alpha \) of \( G_{p} \)
    may be transformed into a \( G_{p} \)-reduced number by 
    finitely many applications of \( TS_{\lambda}^{-r} \), where 
    at each step 
    \( r=\left[\frac{\alpha}{\lambda}\right]+1 \).
    Furthermore, 
    \( TS_{\lambda}^{-r} \) 
    maps reduced numbers to reduced numbers, so 
    the reduced numbers fall into disjoint cycles.
    Finally, every equivalence between \( G_{p} \)-reduced numbers is 
    obtained by iteration of \( TS_{\lambda}^{-r} \);
    thus each equivalence class of hyperbolic 
    fixed points contains one cycle of \( G_{p} \)-reduced numbers.
\end{theorem}

\begin{proof}
    Suppose that \( \alpha = \alpha_{0} \) is a hyperbolic fixed point of \( G_{p} \).
    Then 
    \[ \alpha_{0} 
        = [r_{0};r_{1},\ldots,r_{n-1},\overline{r_{n},\ldots,r_{n+m-1}}], \]
    \( r_{0}=\left[\frac{\alpha_{0}}{\lambda}\right]+1 \).
    Put 
    \begin{eqnarray*}
        \alpha_{1} & = & TS^{-r_{0}}\alpha_{0}  \\
         & = & [r_{1}; r_{2}, \ldots, r_{n-1}, 
                  \overline{r_{n},\ldots,r_{n+m-1}}],
    \end{eqnarray*}
    with \( S = S_{\lambda} \).
    We repeat this process, at each step using the mapping
    \( TS^{-r_{j}} \), 
    \( r_{j}=\left[\frac{\alpha_{j}}{\lambda}\right]+1 \)
    to calculate \( \alpha_{j+1} \).
    After \( n \) steps, we have 
    \begin{eqnarray*}
        \alpha_{n} & = & TS^{-r_{n-1}}\alpha_{n-1}  \\
           & = & [\overline{r_{n};\ldots,r_{n+m-1}}],
    \end{eqnarray*}
    which is reduced.
    
    Next, suppose that \( \beta=\beta_{0} \) is a reduced number, so
    \[ \beta_{0} = [\overline{r_{0};r_{1},\ldots,r_{n}}], \] 
    \( r_{0}=\left[\frac{\beta_{0}}{\lambda}\right]+1 \).
    Then 
    \begin{eqnarray*}
        \beta_{1} & = & TS^{-r_{0}}\beta_{0}  \\
         & = & [\overline{r_{1};r_{2},\ldots,r_{n},r_{0}}],
    \end{eqnarray*}
     and \( \beta_{1} \) is also reduced.
    Repeating this process yields a cycle of \( n \) reduced numbers.
    Since each mapping is uniquely determined, cycles of reduced 
    numbers must be disjoint.
    
    The remaining statements in the Theorem follow from Lemma 
    \ref{lem.EquivFixedPts}.
\end{proof}

The proof of Theorem \ref{thm.ReductionofNumbers} makes it clear that 
there is a \( 1-1 \) correspondence between cycles of reduced forms 
and admissible periods of hyperbolic numbers.
We restate Theorem \ref{thm.ReductionofNumbers} for reduction of hyperbolic 
\( \lambda \)-BQFs.
\begin{corollary}\label{cor.ReductionofForms}
    Let \( p \geq 3 \) and put
    \( \lambda = \lambda_{p} \).
    Every hyperbolic \( \lambda \)-BQF \( Q \) 
    may be transformed into a reduced \( \lambda \)-BQF by 
    finitely many applications of \( S_{\lambda}^{r}T \), where 
    at each step  
    \( r=\left[\frac{\alpha_{Q}}{\lambda}\right]+1 \).
    Furthermore, \( S_{\lambda}^{r}T \) maps reduced forms to reduced 
    forms, so the reduced forms fall into disjoint cycles.
    Finally, every equivalence between reduced \( \lambda \)-BQFs is 
    obtained by iteration of \( S_{\lambda}^{r}T \);
    thus each equivalence class of \( \lambda \)-BQFs with hyperbolic 
    forms contains one cycle of reduced forms.
\end{corollary}

\begin{example}\label{ex.ReducingFormsandNumbers}
    Put \( \lambda = \lambda_{5} = \frac{1+\sqrt{5}}{2} \),
    \( S = S_{\lambda} \), and let
    \( \alpha_{0} = [2;3,\overline{2,1,1,4}] \).
    Then 
    \( M_{\alpha_{0}} 
        = \bigl( \begin{smallmatrix} -9\lambda -6 & \; 51\lambda+32 \\ 
	  -3\lambda -2 & \; 18\lambda +9 \end{smallmatrix} \bigr) \)
    generates the stabilizer of
    \( \alpha_{0} \), and
    \( Q_{\alpha_{0}} 
             = [-3\lambda -2,27\lambda +15,-51\lambda -32]\) 
    is the \( \lambda_{5} \)-BQF corresponding to \( \alpha_{0} \).
    Now \( Q_{\alpha_{0}} \) is in a hyperbolic equivalence class 
    \( \mathcal{A} \) of forms 
    of discriminant \( D = 135\lambda + 86 = \frac{307+135\sqrt{5}}{2} \).
    We reduce \( Q_{\alpha_{0}} \) by 
    \begin{eqnarray*}
        Q_{\alpha_{1}} & = & Q_{\alpha_{0}} \circ S^{2}T  
	            = [\lambda +2,-7\lambda -3,-3\lambda -2],  \\
        Q_{\alpha_{2}} & = & Q_{\alpha_{1}} \circ S^{3}T 
	            = [3\lambda +4,-11\lambda -3,\lambda +2],
    \end{eqnarray*}
    which is reduced.
    The cycle of reduced \( \lambda_{5} \)-BQFs in \( \mathcal{A} \) is
    \begin{eqnarray*}
        Q_{\alpha_{3}} & = & Q_{\alpha_{2}} \circ S^{2}T 
                    = [13\lambda +8,-17\lambda -9,3\lambda +4],  \\
        Q_{\alpha_{4}} & = & Q_{\alpha_{3}} \circ ST 
                    = [11\lambda +8,-25\lambda -17,13\lambda +8],  \\
        Q_{\alpha_{5}} & = & Q_{\alpha_{4}} \circ ST 
                        = [\lambda +2,-13\lambda -5,11\lambda +8], 
                    \text{and}  \\
        Q_{\alpha_{2}} & = & Q_{\alpha_{5}} \circ S^{4}T 
                    = [3\lambda +4,-11\lambda -3,\lambda +2].
    \end{eqnarray*}
    The successive values of \( r \) in \( S^{r}T \) are the successive entries 
    in the \( \lambda_{5} \)-CF for \( \alpha_{0} \).
    The corresponding reduced numbers are
    \begin{eqnarray*}
	\alpha_{2} & = & TS^{-3}\alpha_{1} 
			= TS^{-3}TS^{-2}\alpha_{0}
	                = [\overline{2;1,1,4}], \\
        \alpha_{3} & = & TS^{-2}\alpha_{2} 
			= [\overline{1;1,4,2}],  \\
        \alpha_{4} & = & TS^{-1}\alpha_{3} 
			= [\overline{1;4,2,1}], \text{and}  \\
        \alpha_{5} & = & TS^{-1}\alpha_{4} 
			= [\overline{4;2,1,1}].
    \end{eqnarray*}
\end{example}

\section{SIMPLE FORMS AND NUMBERS}
\label{sec.Simple}
\hspace{\parindent}

In this section we define simple \( \lambda \)-BQFs and 
simple numbers, which are easily characterized and are related to 
reduced forms and numbers.
We use the function \( \Phi_{p} \), defined in Section 
\ref{sec.Heckegroups}, to put the simple forms in each hyperbolic 
equivalence class into a cycle.

We call a hyperbolic \( \lambda_{p} \)-BQF \( Q = [A,B,C] \) 
\( G_{p} \)-\emph{simple} if
\( A > 0 > C \).
If \( Q \) is a simple form, we will say that \( \alpha_{Q} \) is a \( G_{p} 
\)-simple number.

\begin{lemma}\label{lem.SimpleNumbersBQFs}
    Let \( p \geq 3 \) and put
    \( \lambda = \lambda_{p} \).
    Suppose that \( Q = [A,B,C] \) is a hyperbolic \( \lambda \)-BQF 
    associated with \( \alpha = \alpha_{Q} \).
    Then \( Q \) is simple if and only if \( \alpha^{\prime}<0<\alpha \).
\end{lemma}
\begin{proof}
    The proof is an exercise in calculating with inequalities.
\end{proof}

Next we establish the connection between reduced and simple numbers.

\begin{theorem}\label{thm.ReducedandSimpleNumbers}
    Let \( p \geq 3 \) and put
    \( \lambda = \lambda_{p} \).
    Suppose that \( \alpha \) is a \( G_{p} \)-simple number.
    Then \( S_{\lambda}^{n}\alpha \),  
    \( n = -\left[\frac{\alpha^{\prime}}{\lambda}\right] \)
    is a \( G_{p} \)-reduced number.
    The set of \( G_{p} \)-simple numbers is given by 
    \[ \mathcal{Z} = \left\{S_{\lambda}^{-i}\beta \left|\right. 
    \beta \text{ is \( G_{p}- \)reduced, } 
        1 \leq i \leq \left[\frac{\beta}{\lambda} \right]\right\}. \]
\end{theorem}

\begin{proof}
    Put \( S = S_{\lambda} \)
    and let \( \alpha \) be a 
    \( G_{p} \)-simple number.
    Then \( \alpha^{\prime}<0<\alpha \), so 
    \( n = -\left[\frac{\alpha^{\prime}}{\lambda}\right] \geq 1 \)
    and \( S^{n}\alpha > \lambda \).
    Also, \( -n<\frac{\alpha^{\prime}}{\lambda} < 1-n \) implies that 
    \( 0 < S^{n}\alpha^{\prime} < \lambda \).
    Thus \( 0<S^{n}\alpha^{\prime}<\lambda <S^{n}\alpha \), and
    \( \beta = S^{n}\alpha \) is a reduced number
    by Theorem \ref{thm.ReducedNumbers}.

    To prove the second statement, we first let \( \alpha \) be any 
    \( G_{p} \)-simple number, so
    \( \alpha^{\prime} < 0 < \alpha \).
    Then \( \beta = S^{i}\alpha \), 
    \( i = -\left[\frac{\alpha^{\prime}}{\lambda}\right] \),
    is reduced and 
    \( \left[\frac{\beta}{\lambda}\right] = 
        \left[\frac{\alpha}{\lambda}\right] + i \geq i \).
    Now 
    \( \left[\frac{\alpha^{\prime}}{\lambda}\right] \leq -1 \), so
    \( i \geq 1 \),
    and \( \alpha \in \mathcal{Z} \).
    
    To finish proving the second statement, we suppose that 
    \( \alpha \in \mathcal{Z} \).
    Then \( \alpha = S^{-i}\beta \), 
    \( 1 \leq i \leq \left[\frac{\beta}{\lambda}\right] \),
    for some reduced \( \beta \).
    Now \( i \leq \left[\frac{\beta}{\lambda}\right] \) implies that
    \( \alpha = S^{-i}\beta > 0 \).
    The fact that \( i \geq 1 \), along with \( \beta^{\prime} < \lambda \), 
    implies that \( \alpha^{\prime} = S^{-i}\beta^{\prime} < 0 \).
    Thus \( \alpha^{\prime} < 0 < \alpha \) and \( \alpha \) is simple.
\end{proof}

We restate Theorem \ref{thm.ReducedandSimpleNumbers} in terms of hyperbolic 
\( \lambda \)-BQFs.
\begin{corollary}\label{cor.ReducedandSimpleForms}
    Let \( p \geq 3 \) and put
    \( \lambda = \lambda_{p} \).
    Every simple \( \lambda \)-BQF \( Q \) is transformed into a 
    reduced form by a single application of 
    \( S_{\lambda}^{-n} \), with 
    \( n = -\left[\frac{\alpha_{Q}^{\prime}}{\lambda}\right] \).
    The set of simple \( \lambda \)-BQFs is given by 
    \[ \left\{Q \circ S_{\lambda}^{i} \left|\right. 
    Q \text{ is \( G_{p}- \)reduced, } 
            1 \leq i \leq \left[\frac{\beta_{Q}}{\lambda}\right]\right\}. \]
\end{corollary}

Since there are a finite number of reduced \( \lambda \)-BQFs of 
a given discriminant, 
the Corollary implies that there are also a finite 
number of simple forms of a given discriminant.
Also, every reduced form \( Q \) (or number \( \alpha \)) is connected 
by \( S^{j} \) to 
\( \left[\frac{\beta_{Q}}{\lambda}\right] \) simple forms (numbers).
In particular, if \( \beta_{Q} < \lambda \) 
then \( Q \circ S^{j} \) fails to be simple for any \( j \).

Simple numbers (and associated forms) may also be put into cycles, using the 
function \( \Phi_{p} \) defined in Section \ref{sec.Heckegroups}.

\begin{theorem}\label{thm.SimpleCycle}
    Let \( p \geq 3 \) and put
    \( \lambda = \lambda_{p} \).
    The finite orbits of \( \Phi_{p} \) are the set \( \{0\} \) and 
    the sets 
    \[ \mathcal{Z}_{\mathcal{A}} 
       = \left\{\alpha_{Q} \mid Q \in \mathcal{A} \text{, simple}\right\}, \]
    where \( \mathcal{A} \) runs over all hyperbolic \( G_{p} 
    \)-equivalence classes of \( \lambda \)-BQFs.
\end{theorem}

\begin{proof}
    Put \( S = S_{\lambda} \) 
    and \( U = U_{\lambda} \).
    Any finite orbit except \( \{0\} \) has the form 
    \begin{multline*}
        \mathcal{C} =
        \{\alpha_{1}, \alpha_{1}-\lambda, \dots, \alpha_{1}-m_{1}\lambda,
	\alpha_{2}, \alpha_{2}-\lambda, \dots,  \\
	\alpha_{2}-m_{2}\lambda, \dots, 
	\alpha_{s}, \alpha_{s}-\lambda, \dots, \alpha_{s}-m_{s}\lambda \}, 
    \end{multline*}
    for some positive real numbers \( \alpha_{1}, \alpha_{2}, \dots, \alpha_{s} \),
    where \( m_{j}=\left[\frac{\alpha_{j}}{\lambda}\right] \geq 0 \),
    \( 1 \leq j \leq s \).
    We must also have that 
    \[ \alpha_{j+1} = TU^{i_{j}}\left(\alpha_{j}-m_{j}\lambda\right), \]
    \( 1 \leq j \leq s-1 \), as well as
    \[ \alpha_{1} = TU^{i_{s}}\left(\alpha_{s}-m_{s}\lambda\right), \]
    with \( 1 \leq i_{j} \leq p-2 \) for all \( j \).
    Now since \( m_{j}=\left[\frac{\alpha_{j}}{\lambda}\right] \),
    we have 
    \( 0 < \alpha_{j}-m_{j}\lambda < \lambda \), so 
    \( U^{\ell_{j}+1}(0) < \alpha_{j}-m_{j}\lambda < U^{\ell_{j}}(0) \)
    for some \( \ell_{j} \), \( 2 \leq \ell_{j} \leq p-1 \).
    We have eliminated the possibility that 
    \( \alpha_{j}-m_{j}\lambda = U^{\ell_{j}+1}(0) \), 
    since then \( \alpha_{j}-m_{j}\lambda \) would not be part of a 
    finite orbit of \( \Phi_{p} \).
    Thus 
    \begin{eqnarray*}
        \alpha_{j+1} & = & \Phi_{p}(\alpha_{j}-m_{j}\lambda)  \\
         & = & TU^{p-\ell_{j}}S^{-m_{j}}\alpha_{j},
    \end{eqnarray*}
    \( 1 \leq j \leq s-1 \), and 
    \begin{eqnarray*}
        \alpha_{1} & = & \Phi_{p}(\alpha_{s}-m_{s}\lambda)  \\
         & = & TU^{p-\ell_{s}}S^{-m_{s}}\alpha_{s}.
    \end{eqnarray*}
    For each \( j \) we put \( \beta_{j} = S\alpha_{j} \) and 
    \( r_{j} = \left[\frac{\beta_{j}}{\lambda}\right]+1 = m_{j}+2 \geq 2 \).
    Then  
    \begin{eqnarray*}
        \beta_{j+1} & = & S\alpha_{j+1}  \\
         & = & U^{p-\ell_{j}+1}S^{-m_{j}}\alpha_{j}  \\
         & = & U^{p-\ell_{j}+1}S^{-m_{j}-1}\beta_{j},
    \end{eqnarray*}
    for \( j \neq s \), and 
    \[ \beta_{0} = U^{p-\ell_{s}+1}S^{-m_{s}-1}\beta_{s}. \]
    Reversing direction and using \( m_{j} = r_{j}-2 \), 
    we have 
    \[ \beta_{j} = S^{r_{j}}TU^{\ell_{j}-2}\beta_{j+1}, \] 
    \( j \neq s \), 
    and
    \[ \beta_{s} = S^{r_{s}}TU^{\ell_{s}-2}\beta_{0}. \]
    Therefore
    \[ \beta_{j} = 
    [\overline{r_{0};\underbrace{1, \dots,1}_{\ell_{j}-2},r_{j+1},\dots}], 
    \]
    \( j \neq s \), 
    and 
    \[ \beta_{s} = 
    [\overline{r_{s};\underbrace{1, \dots,1}_{\ell_{s}-2},r_{0},\dots}]. \]
    Now \( 0 \leq \ell_{j}-2 \leq p-3 \),
    so the \( \lambda \)-CF expansions are all admissible and \( \beta_{j} \) 
    is \( G_{p} \)-reduced for each \( j \).
    Clearly, the \( \beta_{j} \) are all part of the same cycle of reduced numbers
    in a \( G_{p} \)-equivalence class of numbers.
    We may
    let \( \mathcal{A} \) represent the equivalence 
    class of \( \lambda \)-BQFs containing the corresponding forms.
    If any ones occur in the \( \lambda \)-CF expansions, the 
    cycle contains reduced numbers other than the \( \beta_{j} \).
    But these other reduced numbers are of the form 
    \( \beta = [\overline{1; \dots}] < \lambda \),  
    and are not connected with any \( G_{p} \)-simple numbers.
    Thus all of the simple numbers associated with forms in 
    \( \mathcal{A} \) are of the form \( S^{-i}\beta_{j} \), 
    \( 1 \leq j \leq s \), where \( i \) is a positive integer.
    By Theorem \ref{thm.ReducedandSimpleNumbers}, 
    \begin{eqnarray*}
        \mathcal{Z}_{\mathcal{A}} & = & \left\{S^{-i}\beta_{j} \mid 
	             1\leq i \leq \left[\frac{\beta_{j}}{\lambda}\right]\right\}
		     _{j=1}^{s}  \\
         & = & \left\{S^{-(i-1)}\alpha_{j} \mid 
	             0\leq i-1 \leq m_{i}\right\}_{j=1}^{s}  \\
         & = & \mathcal{C}.
    \end{eqnarray*}
\end{proof}

\begin{example}\label{ex.SimpleFormsandNumbers}
    In Example \ref{ex.ReducingFormsandNumbers} we found the reduced 
    \( \lambda_{5} \)-BQFs in a \( G_{5} \)-equivalence class 
    \( \mathcal{A} \).
    The forms 
    \[ Q_{\alpha_{2}} = [3\lambda +4,-11\lambda -3,\lambda +2], \]
    and 
    \[ Q_{\alpha_{5}} = [\lambda +2,-13\lambda -5,11\lambda +8], \]
    are the only reduced forms 
    \( Q = [A,B,C] \) in 
    \( \mathcal{A} \)
    which correspond to reduced numbers greater than \( \lambda \).
    The simple forms in \( \mathcal{A} \) are all related to \( Q_{\alpha_{2}} \) 
    and \( Q_{\alpha_{5}} \) by \( S^{j} \) as
    \begin{eqnarray*}
        Q_{\alpha_{2}} \circ S & = & [3\lambda +4,3\lambda +3,-3\lambda -2],  \\
        Q_{\alpha_{5}} \circ S & = & [\lambda +2,-7\lambda -3,-3\lambda -2],  \\
        Q_{\alpha_{5}} \circ S^{2} & = & [\lambda +2,-\lambda -1,-9\lambda -6], 
	            \text{and} \\
        Q_{\alpha_{5}} \circ S^{3} & = & [\lambda +2,5\lambda +1,-7\lambda -4].
    \end{eqnarray*}
    We could calculate the corresponding simple numbers directly from 
    these \( \lambda_{5} \)-BQFs.
    Instead, we will use the reduced numbers we found in 
    Example \ref{ex.ReducingFormsandNumbers} along with 
    Theorem \ref{thm.ReducedandSimpleNumbers}.
    The simple numbers in \( \mathcal{A} \) are
    \begin{eqnarray*}
        \mathcal{Z_{A}} 
	 & = & \left\{ S^{-1}\alpha_{2}, S^{-1}\alpha_{5}, 
	            S^{-2}\alpha_{5},
		    S^{-3}\alpha_{5}\right\}  \\
         & = & \left\{ [1;\overline{1,1,4,2}],
	               [3;\overline{2,1,1,4}],
		       [2;\overline{2,1,1,4}],
		       [1;\overline{2,1,1,4}] \right\}  \\
         & = & \left\{ \frac{-3\lambda -3+\sqrt{D}}{6\lambda +8},
	               \frac{7\lambda +3+\sqrt{D}}{2\lambda +4},
		       \frac{\lambda +1+\sqrt{D}}{2\lambda +4},
		       \frac{-5\lambda -1+\sqrt{D}}{2\lambda +4} \right\} \\
	 & = & \left\{ \alpha_{1}, \alpha_{2}, \alpha_{3}, \alpha_{4} \right\}, \\
    \end{eqnarray*}
    where \( D = 135\lambda + 86 = \frac{307+135\sqrt{5}}{2} \).
    Finally, these simple numbers form a finite orbit of 
    \( \Phi_{5} \).
    We have
    \begin{eqnarray*}
        \Phi_{5}(\alpha_{1}) & = & TU\alpha_{1} = \alpha_{2},  \\
        \Phi_{5}(\alpha_{2}) & = & TU^{4}\alpha_{2} = \alpha_{2}-\lambda = \alpha_{3},  \\
        \Phi_{5}(\alpha_{3}) & = & TU^{4}\alpha_{3} = \alpha_{3}-\lambda = \alpha_{4},  \\
        \Phi_{5}(\alpha_{4}) & = & TU^{3}\alpha_{4} = \alpha_{1}.
    \end{eqnarray*}
\end{example}

%\pagebreak

\end{document}